\documentclass[12pt]{amsart}
\usepackage[utf8]{inputenc}
\usepackage[T2A]{fontenc}
\usepackage{graphicx, color, hyperref}%
 \usepackage{latexsym, amsmath, amsfonts,  amssymb, euscript, multicol, bm}
\usepackage{mathrsfs}
\usepackage{float,  graphicx}

\usepackage{environ}

\newtheorem{thm}{\bf Theorem}[section]
\newtheorem{cor}[thm]{\bf Corollary}
\newtheorem{lem}[thm]{\bf Lemma}
\newtheorem{prop}[thm]{\bf Proposition}
\newtheorem{dfn}[thm]{\bf Definition}

\title{Finiteness properties for self-similar continua}
\author{Andrei Tetenov}
\address{Gorno-Altaisk State University and Novosibirsk State University}

\begin{document}
\maketitle

\newcommand{\rr}{\mathbb{R}}
\newcommand \nn {\mathbb{N}}
\newcommand \zz {\mathbb{Z}}
\newcommand \bbc {\mathbb{C}}
\newcommand \rd {\mathbb{R}^d}

 \newcommand {\al} {\alpha}
\newcommand {\be} {\beta}
\newcommand {\da} {\delta}
\newcommand {\Da} {\Delta}
\newcommand {\ga} {\gamma}
\newcommand {\Ga} {\Gamma}
\newcommand {\la} {\lambda}
\newcommand {\La} {\Lambda}
\newcommand{\om}{\omega}
\newcommand{\Om}{\Omega}
\newcommand {\sa} {\sigma}
\newcommand {\Sa} {\Sigma}
\newcommand {\te} {\theta}
\newcommand {\fy} {\varphi}
\newcommand {\ep} {\varepsilon}
\newcommand{\e}{\varepsilon}
\newcommand{\eps}{\epsilon}

\newcommand{\VEC}{\overrightarrow}
\newcommand{\Ra}{\Rightarrow}
\newcommand{\IN}{{\subset}}
\newcommand{\NI}{{\supset}}
\newcommand \dd  {\partial}
\newcommand {\mmm}{{\setminus}}
\newcommand{\probel}{\vspace{. 5cm}}
\newcommand{\8}{{\infty}}
\newcommand{\0}{{\varnothing}}
\newcommand{\vse}{$\blacksquare$}
\newcommand{\ov}{\overline}
\newcommand{\ia}{{I^*}}
\newcommand{\io}{{I^\infty}}
\newcommand{\pim}[1]{{\pi^{-1}(#1)}}
\newcommand{\rmi}{{r_{min}}}
\newcommand{\Id}{{\rm Id}}

\newcommand {\bfep} {{{\bar \varepsilon}}}
\newcommand {\Dl} {\Delta}
\newcommand{\vA}{{\vec {A}}}
\newcommand{\vB}{{\vec {B}}}
\newcommand{\vF}{{\vec {F}}}
\newcommand{\vf}{{\vec {f}}}
\newcommand{\vh}{{\vec {h}}}
\newcommand{\vJ}{{\vec {J}}}
\newcommand{\vK}{{\vec {K}}}
\newcommand{\vP}{{\vec {P}}}
\newcommand{\vX}{{\vec {X}}}
\newcommand{\vY}{{\vec {Y}}}
\newcommand{\vZ}{{\vec {Z}}}
\newcommand{\vx}{{\vec {x}}}
\newcommand{\va}{{\vec {a}}}
\newcommand{\tga}{{\tilde {\gamma}}}
\newcommand{\tda}{{\tilde {\delta}}}

\newcommand{\hf}{{\hat {f}}}
\newcommand{\hg}{{\hat {g}}}

\newcommand{\bj}{{\bm {j}}}
\newcommand{\bi}{{\bm {i}}}
\newcommand{\bk}{{\bm {k}}}
\newcommand{\bl}{{\bm {l}}}
\newcommand{\bma}{{\bm {\alpha}}}
\newcommand{\bmb}{{\bm {\beta}}}
\newcommand{\bX}{{\bm {X}}}

\newcommand{\eA}{{\EuScript A}}
\newcommand{\eB}{{\EuScript B}}
\newcommand{\wP}{{\widetilde P}}
\newcommand{\eU}{{\EuScript U}}
\newcommand{\eM}{{\mathcal M}}
\newcommand{\eV}{{\EuScript V}}
\newcommand{\eS}{{\EuScript S}}
\newcommand{\eH}{{\EuScript H}}
\newcommand{\eC}{{\EuScript C}}
\newcommand{\eP}{{\EuScript P}}
\newcommand{\eT}{{\EuScript T}}
\newcommand{\eG}{{\EuScript G}}
\newcommand{\eK}{{\EuScript K}}
\newcommand{\eQ}{{\EuScript Q}}
\newcommand{\eF}{{\EuScript F}}
\newcommand{\eZ}{{\EuScript Z}}
\newcommand{\eL}{{\EuScript L}}
\newcommand{\eD}{{\EuScript D}}
\newcommand{\E}{{\EuScript E}}
\newcommand{\klm}{\mathfrak{m}}
\def \diam {\mathop{\rm diam}\nolimits}
\def \fix {\mathop{\rm fix}\nolimits}
\def \Lip {\mathop{\rm Lip}\nolimits}
\def \Ord {\mathop{\rm Ord}\nolimits}
\def \Arg {\mathop{\rm Arg}\nolimits}

\newcommand{\beq}{\begin{equation}}
\newcommand{\eeq}{\end{equation}}

\newcommand{\dok}{{\bf{Proof: }}}
\newcommand{\red}{\textcolor{red}}
\newcommand{\yellow}{\textcolor{yellow}}
\newcommand{\blue}{\textcolor{blue}}
\newcommand{\green}{\textcolor{green!60!black}}

\begin{abstract}
 We consider self-similar continua possessing finite intersection property and  
 prove  intersection graph criterion for self-similar dendrites,  finite order Theorem for such  continua  satisfying open set condition in $\rr^n$, and  parameter matching Theorem   which  states that  all Jordan arcs  starting from a intersection point in such continuum on a plane should have the same slope parameter at that point. 
\end{abstract}


\bigskip

\section*{Introduction}

Let $\eS=\{S_1, . . ., S_m\}$ be a system of injective contraction maps in a complete metric space $X$.  A non-empty compact set $K$ satisfying the equation  $K=\bigcup\limits_{i=1}^m S_i(K)$ is called {\em the attractor of the system} $\eS$ and the sets $K_i=S_i(K)$,  where $i=1, . . ., m$,  are called {\em the pieces} of the set $K$.  We will focus mostly on a case in which the space $X$ is $\rr^n$,  the maps $S_i$ are the similarities of $\rr^n$ and the attractor $K$ is connected.  In this case we say $K$ is a {\em self-similar continuum}.

 We say the self-similar continuum $K$ (as well as the system $\eS$) has  finite intersection property  if for any non-equal $i, j$,  the intersection of pieces $K_i\cap K_j$ is finite. If $s$ is an upper bound for all $\#(K_i\cap K_j)$, we say that $\eS$ is a {\em FI(s)-system of contractions}. \\

For a long time it seemed highly likely that finite intersection property could imply open set condition (see Definition \ref{OSC}),  at least in case of one-point intersections.   C.~Bandt and H.~Rao proved in \cite{BR} that FI(s)-systems of similarities in $\rr^2$ with connected attractor satisfy OSC.  From the other side,  it was proved in  \cite{TKV} that in $\rr^3$ this does not hold  and in \cite{nonlin}  it was also shown that there are one-point intersection systems (with totally disconnected attractor) in $\rr$ which violate OSC.  Therefore the interplay between open set condition and finite intersection property is a bit more intricate than it could be supposed.   

 
In this paper we prove three theorems which are  useful tool for analyzing self-similar continua which possess finite intersection pro\-perty. \\

The first is the intersection graph criterion for self-similar dendrites (Theorem \ref{dtest}). It states that if  $\eS$  is a FI(s)-system of contractions in a complete metric space $X$, then its attractor $K$ is a dendrite if and only if the intersection graph $\Ga(\eS)$ of the system $\eS$ is a tree.

It has a  simple topological proof and is applicable to different kinds of systems $\eS=\{S_1, . . ., S_m\}$ of injective contractions in complete metric  spaces. \\

The second one is  finite order Theorem \ref{ddfin} for FI(s)-self-similar continua in $\rr^n$ which  satisfy open set condition. It says there is a  finite uniform upper  bound for the cardinality  of the set of addresses of boundary points,  for the number $N_C$ of components of $V\mmm \{x\}$ (resp.  $V\mmm K_\bj$),  where $V$ is a connected neighborhood of a given set,  and for the topological order of any  point $x\in K$ and of any piece $K_\bj$.

The third  is parameter matching Theorem \ref{parmat} for the intersection points, which extends our result for polygonal dendrites \cite[Theorem 3.12]{Adam} to FI(s)-self-similar continua. It states that if a point $p\in K_i\cap K_j$  has a periodic address, then all Jordan arcs $\ga\IN K$ with the end point $p$ have the same slope parameter at the point $p$. \\

In Section 2 we define  a FI(s)-system $\eS$  of contractions  and   its intersec\-tion graph  $\Ga(\eS)$.  This graph is a bipartite graph whose "white" vertices are the pieces $K_i$ and "black" vertices are the intersection points $p\in K_i\cap K_j$.   In case when $\Ga(\eS)$ is a tree we prove  Theorem \ref{dtest}.

In Section 3 we consider FI(s)-systems  $\eS$  of contracting similarities which satisfy open set condition and prove  Theorem  \ref{ddfin}. 
As an outcome of this Theorem we find that  each $x\in K$ and each $K_\bj\IN K$ has a fundamental system of {\em stable neighborhoods} $V_\bj$ (Theorem \ref{nbh}). Finally we prove a ramification order estimate for self-similar dendrites which need not satisfy finite intersection property (Theorem 2.4).

In Section 4 we consider the fixed points $z_\bj$ of the similarities $S_\bj\in G_\eS$ and prove the existence of invariant arcs in each of the components $Q_k$ of the set $K\mmm \{z_\bj\}$. This yields the
parameter matching theorem for the points $x\in K$  which have multiple preperiodic addresses. 

\newpage

\section{ Preliminaries }

\begin{dfn} 
Let $\eS=\{S_1,  S_2, \ldots,  S_m\}$ be a system of injective contractions on a complete metric space $(X,  d)$. 
 A nonempty compact $K\IN X$ is called the attractor of the system $\eS$,  if $K =\bigcup\limits_{i = 1}^m S_i (K)$. \end{dfn}
 
  The system $\eS$ defines its {\em Hutchinson operator} $T$  by the equation $T(A) =\bigcup\limits_{i = 1}^m S_i (A)$.  By Hutchinson's Theorem \cite{Hut},  the attractor $K$ always exists and is unique for $\eS$.  We also call the   set $K$ {\em self-similar} with respect to $\eS$.  

Given the system $\eS$,  the set $I=\{1, 2, . . ., m\}$ is called the {\em set of indices},  the set $\ia=\bigcup\limits_{n=1}^\8 I^n$ 
is the set of  {\em multiindices} $\bj=j_1j_2. . . j_n$.   
 
If $\bj\in\ia$, we write $S_\bj=S_{j_1}S_{j_2}. . . S_{j_n}$ and denote  $S_\bj(K)$ by $K_\bj$.   
  
The set of all infinite sequences $I^{\8}=\{{\bm\al}=\al_1\al_2\ldots, \al_i\in I\}$ is called {\em the
index space} of the system $\eS$; the map $\pi:I^{\8}\rightarrow K$ which   sends each $\bma$ to  the point $\pi(\bma)=\bigcap\limits_{n=1}^\8 K_{\al_1\ldots\al_n}$ is called the {\em index map}.  If $\pi(\bma)=x$,  then $\bma$ is called {\em an address} of the point $x$.  \\

By $\bi\bj$ we denote the concatenation of respective multi-indices.  We say $\bi\sqsubset\bj$,  if  $\bj=\bi\bl$ for some $\bl\in\ia$; if $\bi\not\sqsubset\bj$ and $\bj\not\sqsubset\bi$,  we say that $\bi$ and $\bj$ are {\em incomparable}. \\
Similarly,  by $\bi\bma$ we denote the concatenation of $\bi$ and $\bma$ and  we write $\bi\sqsubset\bma$,  if  $\bma=\bi\bmb$ for some $\bmb\in\ia$. If $\bma=\bi\bma$, we write $\bma=\bar\bi$ and say that $\bma$ is periodic. We say that $\bma$ is preperiodic if $\bma=\bj\bar\bi$ for some $\bi,\bj\in \ia$.

We denote $\bi\io=\{\bi\bma, \bma\in \io\}$,  and $J\io=\{\bi\bma:\bi\in J, \bma\in \io\}$
 Along with the system $\eS$ we consider its {\em n-th refinement} $\eS^{(n)}=\{S_\bj, \bj\in I^n\}$.   The Hutchinson operator of the system $\eS^{(n)}$ is equal to $T^n$.

 \begin{dfn}\label{jleq}
Let $J'\IN I^*$ be a set of incomparable multiindices,  and let $A\IN I^\8$ be a set of addresses.  We write $J\prec J'$ (resp.  $J\prec A$),  if  $J$ is a set of incomparable multiindices and if there is a surjection $\psi:J'\to J$ or $\psi:A\to J$ such that  $\bj'\sqsupset\psi(\bj')$ (resp.  $\al\sqsupset\psi(\al)$). \\
\end{dfn}

Suppose $B\IN K$ and $J\prec\pim B$. 
 We denote by $V_J(B)$ the set $\bigcup\limits_{\bj\in J}K_\bj$. 

\begin{lem}\label{VJ} If $B\IN K$ is a compact set and $J\prec\pim B$,  then $V_J$ is a neighborhood of $B$ in $K$. \end{lem}
{\bf Proof. } Notice that the family $\{\bj\io,  \bj\in J\}$ is an open covering of $\pim B$.  Since the multiindices $\bj\in J$ are incomparable,  this covering is a disjoint one and therefore the set $J$ is finite.  There is a finite $J'\IN \ia$ such that $J'\io= \io\mmm  J \io$.  The set $\pi(J'\io)$ is compact and $\pi(J'\io)\cap B=\0$.  There is such $\e>0$ that for any $x\in B$ and $y\in \pi(J'\io)$,  $d(x, y)>\e$.  Therefore $V_J(B)$ contains a $\e$-neighborhood of $B$ in $K$,  therefore $V_J(B)$ is a neighborhood of $B$. \vse

 \begin{dfn}
Let $S$ be a contracting injective map on a metric space $X$ and $x=\fix(S)$. 
A Jordan arc $\ga\IN X$ with an endpoint  at $x$ is called an  invariant arc for $S$,  if there is $k\in\nn$ such that  $S^k(\ga)\IN \ga$. 
\end{dfn}

\begin{dfn}\label{OSC} The system ${\eS}$ satisfies the {\em open set condition} (OSC) if there exists a non-empty open set $O\IN X$ such that the sets $S_i (O), \{1\le i\le m\}$ are pairwise disjoint and are contained in $O$. \end{dfn}
For any $\bi, \bj\in \ia$,  $\bi\sqsubset\bj$ iff $S_\bi(O)\NI S_\bj(O)$ and  $\bi$ and $\bj$ are incomparable,  iff $S_\bi(O)\cap S_\bj(O)=\0$. 
The union $\eC$ of all intersections $S_i(K)\cap S_j(K)$,  $i, j\in I,  i\neq j$ is called  {\em the critical set} of the system $\eS$.  The set of all predecessors of the points in $\eC$, 
 $\dd K=\{x\in K: \mbox{ for some }\bj\in\ia, \ \  S_\bj\in \eC$ is called the {\em self-similar boundary} of the set $K$.

 We denote  by $G_\eS=\{S_\bj,  \bj\in\ia\}$ the semigroup,  generated by $\eS$.  In case when $\eS$ is a system of similarities in $\rr^n$ then $\eF=G_\eS^{-1}\circ G_\eS$,  or a set of all compositions $S_{\bf j}^{-1}S_{\bf i}$,  ${\bf i},  {\bf j} \in \ia$,  is the {\em associated family of similarities}\cite{SSS7}.  The system $\eS$ has the {\em weak separation property} (WSP) iff $\rm Id \notin \overline{\eF \setminus \rm Id}$\cite{Lau}, \cite{Zer}. \\

{\bf Definition of $\eM_a$. } According to Zerner's Theorem \cite{Zer},  if the system  $\eS$ of contraction similarities satisfies the Weak Separation Condition,  then for any $a>0$ there is a positive number
\beq \label{m1} \eM_a=\sup\limits_{U\IN\rr^n}\#\{S_\bj:a|U|\rmi <|K_\bj|\le a|U|\ \ \&\ \ \  K_\bj\cap U\neq\0\}. \eeq
If the system $\eS$ possesses the finite intersection property it has no exact overlaps,  and $S_\bi=S_\bj$ iff $\bi=\bj$,  then \eqref{m1} becomes
\beq\label{m2} \eM_a=\sup\limits_{U\IN\rr^n}\#\{\bj:a|U|\rmi <|K_\bj|\le a|U|\ \ \&\ \ \  K_\bj\cap U\neq\0\}. \eeq

\section{The  intersection graph  of a system with finite intersection property} 

In this section we start with  a definition of a FI(s)-set system $\eK$ and define its intersection graph $\Ga(\eK)$. Then we  define the $(\eL, f)$-refinement operation for such systems and  show that the class of tree-like systems is closed under this operation. After all,  we apply the established  properties to systems $\eS$ of  contractions to prove  Theorem  \ref{dtest}. 
\begin{dfn}\label{fipss}
Let $\eK=\{K_i, i\in I=\{1, \ldots, m\}\}$ be а  finite system of continua  in a Hausdorff topological space $X$.  We say $\eK$ possesses {\em finite intersection  property},  if for any  $i\neq j\in I $,  the intersection $P_{ij}=K_i\cap K_j$ is finite.  
 We call $\eK$  a {\em FI(s)-set system},  if $\max\#P_{ij}=s$. 
\end{dfn}
In the settings of the Definition  \ref{fipss} we denote $K=\bigcup\limits_{i\in I}K_i$,   $P=\bigcup\limits_{i\neq j}P_{ij} $ and $P_i=\bigcup\limits_{j\in I\mmm\{i\}}P_{ij}$.  Considering $K$ as the subspace of $ X$ provided by the induced topology,  we see that  the set $P_i$ is the boundary $\dd K_\bi$ of the set $K_i$ in $K$,  and that its interior is $\dot K_i=K_i\mmm P_i$.  Observe that for any $i\in I$,  $\#\dd K_i\le (m-1)s$. 
 \begin{dfn}\label{fipcs}
Let  $\eS=\{S_1, . . ., S_m\}$ be a system of injective contraction maps  on a complete metric space $X$ and $K$ be its attractor.  Let $\eK(\eS)=\{K_1, . . .,  K_m\}$.  
   $\eS$ is called a FI(s)-system of contractions if  the system $\eK(\eS)$ is a FI(s)-set system. \end{dfn}
   
   Along with the FI(s)-system $\eS$ we consider its $n$-th refinement $\eS^n=\{S_\bi, \bi\in I^n\}$ and the system $\eK_n(\eS)=\{K_\bi:\bi\in I^n\}$.  For each $K_\bi$ its  boundary is
   $\dd K_\bi= \bigcup\limits_{\bj\in I^n\mmm \{\bi\}}K_\bj\cup K_\bi$.  We define $P^n=\bigcup\limits_{\bi\in I^n}\dd K_\bi $,  and  the set of all boundary points of the pieces of $K$ is $P^*=\bigcup\limits_{\bi\in I^*}\dd K_\bi$. 
   
   If $\bi, \bj\in \ia$ and $i_1\neq j_1$   then $\#(K_\bi\cap K_\bj)\le \#(K_{i_1}\cap K_{j_1})\le s$.   This implies that  for any incomparable $\bi, \bj\in \ia$   we obtain $\#(K_\bi\cap K_\bj)\le s$.  Therefore   for each $n$   the system $\eS^n$ is also  a FI(s)-system of contractions.

{2. } For a FI(s)-set system $\eK$ (resp.  FI(s)-system $\eS$) we define its {\em intersection graph}  $\Ga(\eK)$ (resp.  $\Ga(\eS)$) as a bipartite graph  $(\eK,  P; E)$ with parts $\eK$ and $P$,  for which an edge  $\{K_i, p\}\in E$ iff  $p\in K_i$.  \\
\begin{figure}[H]
\includegraphics[width=0.5\textwidth]{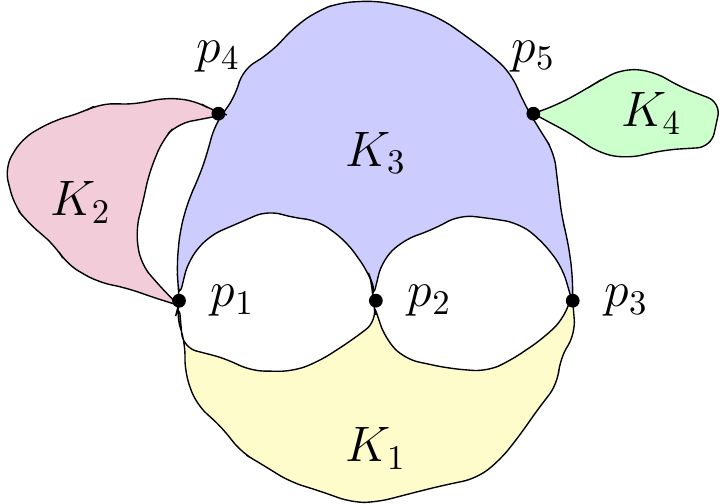}\qquad
\includegraphics[width=0.4\textwidth]{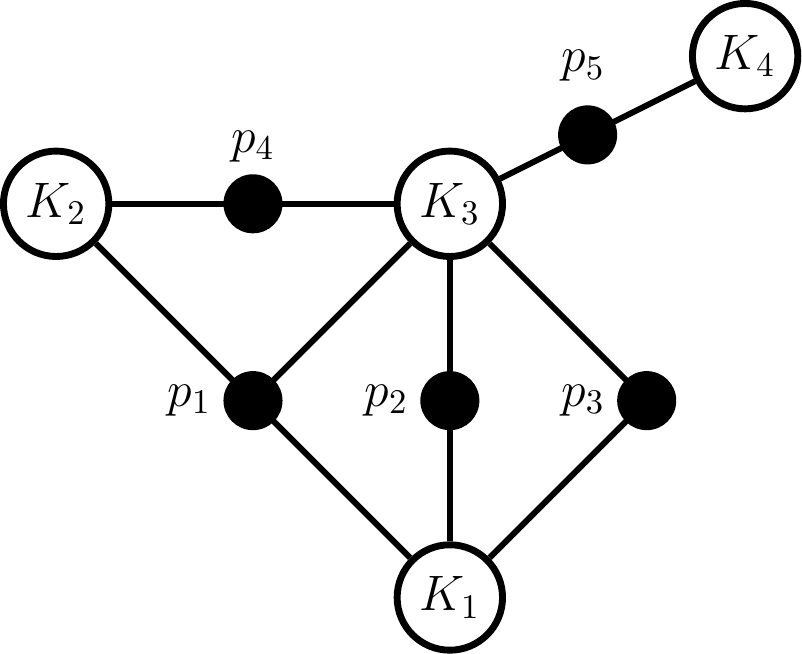} \caption{A FI(3)-set system and its intersection graph.}
\end{figure}

 We call $K_i\in \eK$ {\em white vertices}  and $p\in P$ -- {\em black vertices} of the graph $\Ga$.   The set $N(K_i)$  of the neighbors of any white vertex $K_i$ is  $P_i$,  whereas for any black vertex $p$,   $N(p)=\{K_i:p\in K_i\}$.  Each $p\in P$ is the intersection point of at least two of the sets $K_i$,  therefore $ \mathrm{deg}(p)\ge 2$. \\
  
 Given two FI(s)-set systems  $\eK=\{K_i, i\in I\}$,  $\eL=\{L_i, i\in J\}$,  we define the following  operation which  produces a $(\eL, f)$-refinement $\eK'$ of $\eK$ whose intersection graph $\Ga(\eK')$ is an expansion of the graph $\Ga(\eK)$:\\
We take some $l\in I$ and some  continuous injection $f:L\to K_l$ of the set 
$L=\bigcup\limits_{\eL}L_j$ which satisfies $P_l\IN f(L)$.  The map $f$  transforms the system $\eL$ to  a FI(s)-set system $\eL_f= \{f(L_i),  i\in J\}$. 
Then the system $\eK'=(\eK\mmm\{K_l\})\cup  \eL_f$ is a FI(s)-set system,  whose intersection graph $\Ga(\eK')$ is an expansion of $\Ga(\eK)$. \\

  { 3. } If the intersection graph $\Ga(\eK)$ of a FI(s)-set system $\eS$ is a tree,  then $s=1$.  Moreover,  if the intersection graphs of FI(1)-set systems  $\eK,  \eL$ are trees,  the same is true  for any  $(\eL, f)$-refinement $\eK'$ of $\eK$:

\begin{prop}\label{comp1}
Let $\eK=\{K_i, i\in I\}$,  $\eL=\{L_i, i\in J\}$,  and $l\in I$.  Let the intersection graphs $\Ga(\eK),  \Ga(\eL)$ be trees. 
Then the intersection graph of any $(\eL, f)$-refinement
$\eK'$ of the system $\eK$ is a tree. 
\end{prop}

 \dok Let $l\in I$ and let $f:L\to K_l$ be a continuous injection which maps the set  such that $P_l\IN f(L)$.

 Consider the graph  $\Ga(\eK\mmm\{K_l\})$ obtained by deletion of the star of $K_l$ in $\Ga(\eK)$.   It is a disjoint union of finite number of connected components $Q_p$ each being a non-degenerate tree containing one of the vertices $p\in P_l$.  Since $f:L\to K_l$ is injective,  there are exactly $n_l=\#P_l$ points $p'=f^{-1}(p),  p\in P_l$,  contained in $\bigcup\limits_{i\in J}L_i$. 

 These points are of two kinds.  First are the points $p'\in f^{-1}(P_k)\cap P_\eL$ which correspond to black vertices of $\Ga(\eL)$.   The second type are the points
$p'\in f^{-1}(P_k)\mmm P_\eL$.  Each of these points is contained in some unique $L_i\in \eL$.  We construct an extension  $\tilde\Ga$ of the  graph $\Ga_\eL$,  adding the points of second kind to the set of black vertices  $P_\eL$ and respective edges $(L_i, p')$ to the edge set $E_\eL$.  Such extension does not produce cycles and the graph $\tilde\Ga$ is again a tree with two parts. 

 Now we identify each of the points $f^{-1}(p_i),  p_i\in P_l$,  with  the point $p_i\in Q_i$.  Thus we paste the tree $Q_i$ to the graph $\tilde\Ga$.  

 The resulting graph $\Ga(\eK')$ is a tree with two parts, $\eK'=\eK\mmm\{K_l\}\cup \eL$
and  $P_{\eK'}=P_{\eK_l}\cup(P_\eL\mmm f^{-1}(P_l))$. 
The degree of each vertex $p\in P_{\eL'}$ is $\ge 2$. \vse

 { 4. }
  The refinement operation can be applied to a FI(s)-system $\eS$ of contractions and to its $n$-th intersection graphs $\Ga_n(\eS)=\Ga(\eK_n(\eS))$.  

\begin{lem}
Let $\eS=\{S_i,  i\in I\}$ be a system of injective contraction maps,  such that the intersection graph  $\Ga_1(\eS)$ is a tree.   Then for any  $n\in \nn$,  the intersection graph $\Ga_n(\eS)$ is a tree. 
\end{lem}

  \dok  Suppose that the intersection graph
 $\Ga_{n-1}(\eS)$  is a tree.  Applying $(\eK, S_\bi)$-refinement operation to the piece $K_\bi$ successively for all $\bi\in I^{n-1}$,  we finally obtain that the intersection graph   $\Ga_n$ is also a tree. \vse
 
If the intersection graph of a FI(1)-set system $\eK$ is a tree,  then a simple loop in $K$ cannot pass through any of the boundary points between the sets $K_i$. So, we get the following result.
\begin{lem} Let $\eK$ be a FI(1)-system of sets and
 $\Ga(\eA)$ be a tree.  Let $\ga$ be a simple closed curve in $K$.  Then there is unique $K_l\in\eK$ such that $\ga\IN K_l$. 
\end{lem}
 
 \dok  Let $p$ be some point in $P$ and let $Q_i, Q_j$ be the components of $K\mmm\{p\}$.  Suppose $\ga$ is a closed curve containing some $a\in Q_i$ and $b\in Q_j$.  Since each path connecting $a$ and $b$ passes through $p$,  the point $p$ is a multiple point of $\ga$.  Therefore if $\ga\cap\dot K_i\neq\0$,  then $\ga\in K_i$. \vse

\begin{thm}\label{dtest}
Let $\eS$ be a system of injective contraction maps in a complete metric space $X$   which possesses finite intersection property.     The attractor   $K$ of the system $\eS$ is a dendrite if and only if the intersection graph of the system $\eS$ is a tree. 
\end{thm}

\dok Let $\ga\in K$ be a simple closed curve.  Since for any $n\in \nn$ the graph $\Ga_n$ is a tree,  there is unique $\bj\in I^n$ such that $\ga\in K_\bj$.  Therefore $|\ga|=0$.  \vse

\section{Finiteness theorems for  ramification points and pieces }

Throughout this and further sections,  unless otherwise is stated,  we restrict ourself to the case when  $\eS=\{S_1, . . ., S_m\}$ is a $FI(s)$-system of contracting similarities in $\rr^d$.  Therefore its attractor $K$ is a continuum and for any non-equal $i, j\in I$,  $\#(K_i\cap K_j)\le s$. \\

Consider the following notation: 

${\#\pi^{-1}(x)}$ (resp.  ${\#\pi^{-1}(\dd K_\bj)}$) denotes the cardinality of the set of addresses of a point $x\in K$ (resp.  of boundary points of a piece $K_\bj$). 

${N_C(x)}$ (resp.  ${N_C(K_\bj)}$) denotes the upper bound for the cardinality of the set of connected components of  $V(A)\mmm A$.  Here $V(A)$ denotes a compact connected neighborhood of $A$,  and $A$ stands for $\{x\}$ or $K_\bj$.  

By ${ \Ord(A, X)}$ we denote the topological order of the set $A$ in the space $X$.  \cite{Kur}, \cite[\S 9. 3]{Nad}\\

We prove that if the FI(s)-system $\eS$ satisfies OSC,  then each of the above numbers has a uniform upper bound independent of the choice of  $K_\bj, \bj\in \ia$ or $ x\in K$:

\begin{thm}\label{ddfin}
 If a FI(s)-system $\eS=\{S_1, . . ., S_m\}$ of similarities  satisfies open set condition,  then   each of the numbers:\\ 
 $\#\pi^{-1}(x)$;\   $\#\pi^{-1}(\dd K_\bj)$;\ \   $N_C(\{x\})$;
 \    $N_C(K_\bj)$;\ \  $\Ord(x, K)$; 
  \   $\Ord(K_\bj, K)$ \\ has a uniform finite upper bound independent of  $x\in K$ and $\bj\in \ia$. 
   \end{thm}

 \dok 
 
(i) Cases $\#\pi^{-1}(x)$ and  $\#\pi^{-1}(\dd K_\bj)$.\\ If $x\notin P^*$,  then $\#\pi^{-1}(x)=1$.  Let $x\in P^*$.  Take some  $\rho>0$ and consider the set \[C_\rho=\{\bj\in I^*: \rho\rmi<|K_\bj|\le\rho\, \ \ K_\bj\cap B(x, \rho/2)\neq\0 \}.\] By \cite[Theorem 1]{Zer},  $\#C_\rho\le \eM_1$.  For any $\rho>0$ and any $\al\in\pim{x}$ there is $\bj\sqsubset\al$ such that $\bj\in C_{\rho}$,  therefore $$\#\pim{x}\le\sup\{\#C_\rho, \rho>0\}\le \eM_1$$
 For each $x\in  \dd K_\bi$ choose  $\bj\in I^*$ such that $x\in K_\bj\cap  K_\bi$  and $|K_\bi|\rmi<|K_\bj|\le |K_\bi|$.  The number of all such $\bj$ is at most $\eM_1$,  whereas the number $\#(K_\bi\cap K_\bj)\le s$.   This implies $\#\dd K_\bi\le \eM_1\cdot s$.  Therefore $\#\pi^{-1}(\dd K_\bi)\le \eM_1^2s$.\\ 
 
 (ii) Cases $N_C(\{x\})$ and   $N_C(K_\bj)$.\\
  Let $\eQ=\{Q_1, . . ., Q_n\}$ be some finite set of connected components of  $W\mmm\{x\}$.  Take such $\rho$,  that for any $Q_k\in\eQ$,  $Q_k\mmm B(x, \rho)\neq\0$.  Each component $Q_k\in\eQ$ contains such $y_k$,  that $d(x, y_k)=3\rho/4$.  Let $\bj_k$ be such that $y_k\in K_{\bj_k}$,  and $\rho\rmi/2<|K_{\bj_k}|\le\rho/2$.  Since $x\notin K_{\bj_k}$, \quad
$K_{\bj_k}\IN Q_k$.   Therefore all $\bj_k$ are incomparable and the number of such $\bj$ is no greater than $\eM_{1/3}$.  Thus, $\#\eQ\le \eM_{1/3}$,  which gives the bound for $N_C(\{x\})$. 

 Notice that $N_C(K_\bi)$ is no greater the sum of all $N_C(x)$,  where $x\in\dd K_\bi$,  we obtain $N_C(K_\bi)\le \eM_{1/3}\eM_1s$. \\

(iii) Cases $\Ord(x, K)$ and   $\Ord(K_\bj, K)$.\\
By Lemma \ref{VJ},  for each $x\in K$ and $J\prec\pi^{-1}(x)$, 
 the set $V_J(x)=\bigcup\limits_{\bj\in J} K_{\bj}$ is a connected neighborhood of the point $x$. 

The same is true for each $\bi\in \ia$.  If $J\prec \pi^{-1}(\dd K_\bi)$,  the set
 $V_J(K_\bi)=\bigcup\limits_{\bj\in J} K_{\bj}\cup K_\bi$ is a connected neighborhood of $K_\bi$. 

 Moreover,  the family $\eV(x):=\{V_J(x):\  J\prec \pim{x}\}$  (resp.  $\eV(K_\bi):=\{V_J(K_\bi):\  J\prec \pim{\dd K_\bi}\}$
  is a neighborhood basis in $K$ for the point  $x$ (resp.  the set $K_\bi$),  which consists entirely of connected sets.  
  
Taking into account the inclusion $\dd V_J\IN \bigcup\limits_{\bj\in J}\dd K_\bj$,  we obtain the inequalities  $\Ord(x, K)\le \eM_1^2s$ and    $\Ord(K_\bi, K)\le \eM_1^3s^2$.  \vse

\begin{prop}\label{nbh}
 Let $\eS=\{S_1, . . ., S_m\}$  be  a FI(s)-system   satisfying OSC. Then
 
   1) for any $x\in K$ there is a family $J\prec\pim{x} $ such that
 for any non-equal $\bj, \bk\in J$,  $K_\bj\cap K_\bk=\{x\}$ and
  the number of connected components of 
$V_J\mmm\{x\}$ is equal to $N_C(x)$;

2) for any $\bi\in \ia$ there is a family $J\prec\pim{\dd K_\bi} $ such that
 for any  $\bj\in J$,  $K_\bj\cap \dd K_\bi$ is a unique point,  the sets $K_\bj\mmm K_\bi$ are disjoint and
 the number of connected components of 
$V_J(x)\mmm K_\bi$ is equal to $N_C(K_\bi)$.  \end{prop} 
\dok 
 Let $ \pim{x}=\{\al_1, . . ., \al_n\}$.  There is $J\prec\pim{x}$  which consists of incomparable $\bj\in J$.  Therefore the set $P_J=\bigcup\limits_{\bj\in J} \dd K_{\bj}$ is finite.  Take  $\rho>0$ such that  $B(x, \rho)\cap P_J=\{x\}$.  If $J'\succ J$ satisfies the relation $J\prec J'\prec \pim{x}$ and for any $\bj\in J'$,  $\diam(K_\bj)<\rho$,  then for any non-equal $\bj, \bk\in J'$,  $K_\bj\cap K_\bk=\{x\}$.  This proves (i). \\

Take $\al\in\pim{x}$. Denote the initial substring of length $k$ in $\al$  by $\bj_k$.  
 
The sequence 
$\{n_k:=\#\eQ_{\bj_k}\}$ is non-decreasing.  By  Theorem \ref{ddfin}, it has an upper bound.   Therefore there is
$k_0$ such that if $k\ge k_0$, then  $n_k=n_{k_0}$.  For any $l>k\ge k_0$,   $Q\in\eQ_{\bj_{k}}$ implies $Q\cap K_{\bj_{l}}\in \eQ_{\bj_{l}}$.  Choosing respective $\bj\sqsubset\al$ for each  $\al\in\pim{x}$ we get the desired $J\prec\pim{x}$. \\

 Let $A=\pim{\dd K_\bi}$.  Since $A$ is finite,  we can take some set of incomparable multiindices $J\prec A$.  Let
$\rho=1/2\min\{d(x, y),  x, y\in P_J\}$.  Take such $J'\prec A$,  that $J\prec J'$ and for any $\bj\in J'$,  $|K_\bj|<\rho$.  Then for any $\bj\in J'$,  the intersection $K_\bj\cap P_J$ is a singleton and therefore is contained in $\dd K_\bi$.  \vse\\

\begin{dfn}
The set $\bigcup\limits_{\bj\in J}\eQ_\bj$ defined by  Proposition \ref{nbh}(ii)   is called a {\em stable set of components} for the point $x$. The set $V_J(x)$ is called a {\em stable neighborhood} of the point $x$. \end{dfn}

If  $\eS$  satisfies the assumptions of   Theorem  \ref{ddfin} and its attractor $K$ is a  dendrite,   then  for any $x\in K$ and any stable neighborhood  
$V_J(x)$ of the point $x$,  $\#\dd V_J=N_C(x)$,  therefore by Theorem \ref{ddfin}(3), 
 $\Ord(x, K)\le \eM_{1/3}$.  \\
 
 Surprisingly,  we can prove a similar estimate without the requirement of finite intersection property: 
\begin {thm}\label{wspram}
Suppose a system  $\eS=\{S_1, . . ., S_m\}$ of contracting similarities in $\rr^n$ satisfies WSP,  and its attractor $K$ is a dendrite. 
Then for any $x\in K$,  $\Ord(x, K)\le \eM_{1/2}$. 
\end{thm}

\dok Let $Q_1, . . ., Q_n$ be some finite set of connected components of $K\mmm\{x\}$.  Let $\rho<\min\limits_{1\le k\le n}\diam(Q_k)$.  For each $1\le k\le n$ take some $z_k\in \dd B(x, \rho)\cap Q_k$.  Take such $\bj_k\in I^*$ that $z_k\in K_{\bj_k}$ and $\diam( K_{\bj_k})<\rho$.  Since $K$   is a dendrite and  $x\notin K_{\bj_k} $,  the sets $K_{\bj_k}$ are disjoint.  Therefore by \cite[Theorem 1]{Zer},   $n\le \eM_{1/2}$ for any $x\in K$. \vse\\

%
\section{Boundary points and parameter matching for FI(s) continua  in the plane }

If a FI(s)-system of similarities $\eS$ satisfies OSC, then the fixed points of maps $S_\bi, \bi\in \ia$ possess some special properties.
 \begin{prop}\label{cycl0}
If a FI(s)-system  $\eS=\{S_1, . . ., S_m\}$ of  similarities in $\rr^d$    satisfies OSC,  and $x=\fix{S_\bi}$ for some $\bi\in\ia$,   then:\\
(i) $\pim{x}=\{\overline{\bi}\}$;\\
(ii) The set $\eQ$ of connected components  of the set $K\mmm\{x\}$ is a stable set of components for the point $x$;\\
(iii)  for any component $Q\in\eQ$ there is a point $y\in S_\bi^{-1}(\dd K_\bi)\cap Q$ and an invariant Jordan arc $\ga\IN Q$ for $S_\bi$  with endpoints   $x$ and $y$.
\end{prop}
\dok (i) Suppose there is  $\al\in \pim{x}$ such that  $\al\neq \overline{\bi}$.  Then for any $k\in \nn$,  $\bi^k\al\in\pim{x}$.  Since all these addresses are different, this contradicts Theorem \ref{ddfin}.   \\
(ii) Let $\eQ=\{Q_k, k=1,..., s\}$ be a stable set of components for the point $x$. The equality
$S_\bi(Q_k)=S_\bi(K)\cap Q_{\sa(k)}$ defines a permutation $\sa$ of the set $\{1,...,s\}$.
 Take such $l$  that $\sa^l=\Id$. Then for any $k=1,...,s$, $S^l_\bi(Q_k)=S^l_\bi(K)\cap Q_k$.
 
  Since $K\mmm\{x\}=\bigcup\limits_{k=1}^s(S_\bi^{-l}(Q_k\cap K))$,  the family $\{S_\bi^{-l}(Q_k\cap K)\}$ is also a stable set of components for  $x$.

(iii) To make the notation more simple we use $\bj=\bi^l$ in our further argument.

Let $Q$ be a component of $K\mmm\{x\}$ and let $D_Q=S_\bi^{-1}(\dd S_\bi(K)\cap Q)$. \\ Then    it follows from (ii) that for any $k\in\nn$, $S_\bi^{k}(D)=\dd S_\bi^{k}(K)\cap Q$. \\

We define a map $\fy: D\to D$ the following way.

Since $Q $ is  arcwise connected,  for any $z\in D$ there is  a Jordan arc $\da\in Q$ with endpoints $z$ and $x$. 
Let $\bar\da(z)$ be the closure of a maximal subarc  of $\da\cap (Q\mmm S_\bj(Q))$ containing $z$. We denote the second endpoint of $\bar\da(z)$ by $z'$ and  define a map $\fy:D\to D'$ by the equation  $\fy(z)=S_\bj^{-1}(z')$.  \\

There are  $n\le\#D$ and $y\in D$ such that $\fy^n(y)=y$.  For any  non-negative integer $k<n$ and $y_k=\fy^{k}(y)$,  $\bar\da(y_k)$ is a Jordan arc in $\overline{Q\mmm S_\bi(Q)}$  connecting $y_k$ and $S_\bj( y_{k+1})$. 

Therefore 
$\ga'=\bigcup\limits_{k=0}^n S_\bj^k(\bar\da(y_k)$ is a Jordan arc connecting $y$ and $S_\bj^n(y)$ which lies in $\overline{Q\mmm  S_\bj^n(Q)}$ and   $\ga=\bigcup\limits_{k=0}^\8 S_\bj^{kn}(\ga')\cup\{x\}$ is the desired invariant arc in $Q$ with endpoints $x$ and $y$.  \vse\\

{\bf Remark. } 
The proof of Proposition \ref{cycl0} is purely combinatorial and is based  on finite intersection property and on the assumption that  both $\pim{x}$ and the set of components of $K\mmm\{x\}$ are finite.  This allows us to extend this Proposition  to  FI(s)-systems of contractions in metric spaces. 
\begin{cor}
If  Let  $\eS=\{S_1, . . ., S_m\}$ be a FI(s)-system of contractions in a complete metric space $X$.  If for any $x\in K$ the set $\pim{x}$ and the set of components of $K\mmm\{x\}$ are finite,  then the statements (i), (ii), (iii) of Proposition \ref{cycl0} are fulfilled. 
\end{cor}

If $\eS$ is a FI(s)-system of similarities on a plane,  then it follows from \cite{BR},  that it satisfies open set condition. In dimension 2, the   invariant arcs are more tame and we show that they obey certain restrictions.   To analyze their behavior,  we use complex notation and consider FI(s)-systems $\eS$ as acting in $\bbc$. 

 Let $z_0\in \bbc$ and let $\ga\IN \bbc\mmm\{z_0\}$ be a Jordan arc  with endpoints $z_1, z_2$.  We denote by $\Da(\Arg (z-a))|_\ga$ the increment of  $\Arg(z-a)$ along the arc $\ga$ as $z$ travels from $z_1$ to $z_2$. \\

\begin{dfn} We say that a Jordan arc $\ga=\ga(z_0,b)$ in $\bbc$  {\em has  a  slope parameter $\la$ at the endpoint $z_0$} if  there is $M>0$ such that for any subarc $\ga'=\ga(z_1,z_2)\IN\dot\ga$ 
the increment $\Da\Arg (z-z_0)|_{\ga'}$ of the argument of $z-z_0$ along the arc $\ga'$ satisfies the inequality
$$|\Da\Arg (z-z_0)|_{\ga'} -\la(\log |z_2-z_0| - \log |z_1-z_0|)|\leq M.$$\end{dfn}
\begin{prop}\label{invarcs} If $\ga(z_0,y)$ is an invariant arc for  $S_\bi$ and $S_\bi^n(\ga)\IN \ga$, then the slope parameter of $\ga$ at the point $z_0$ is equal to
\[\la=\dfrac{\Da(\Arg (z-z_0))|_{\ga\mmm S^n_\bi(\ga)}}{k\log\Lip{S_\bi}}.\]
If the point $z_0$ is contained in  boundary $\dd A$  of some unbounded component $A$ of $\bbc\mmm K$, then for any arc $\ga'\IN K$ 
 whose  endpoint is $z_0$, this arc has the slope parameter $\la$ at  $z_0$.
\end{prop}

\dok Replacing, if necessary, the  multiindex $\bi$ by its  n-th multiple, we suppose from now on that $S_\bj(\ga)\IN \ga$ whenever we deal with fixed points and their invariant arcs. 

Denote $R=-\log\Lip{S_\bi}$ and $\al=\Da(\Arg (z-z_0))|_{\ga\mmm S_\bi(\ga)}$ and  let $T$ be a translation of the  complex plane defined  by $T(w)=w-R-i\al$. \\

We slightly modify the proof of \cite[Lemma 3.1]{ATK}.  Consider the map $\fy$ defined by the formula $z=z_0+e^w$ of the complex plane $w=\rho+{i\theta}$ as a universal covering map of $\bbc\mmm \{z_0\}$. 
For any  lifting  $\tga$  of the arc $\ga$ with respect to $\fy$,  $T(\tga)\IN\tga$. Fix one of the liftings and let $\tilde y=\rho_0+i\te_0\in \fy^{-1}(y)$ be its endpoint and let  $\tda=\overline{\tga\mmm T(\tga)}$ be the respective lifting of $\da$. 
There is  $M>0$ such that the arc $\tda$ lies in a strip $|\te-\te_0-\la(\rho-\rho_0)|<M$. By its $T$-invariance, the same is true for $\tga$, thus proving the first statement of the Proposition.

 Now let $\ga'$ be a Jordan arc in $A$, whose endpoints are  $z_0$ and some $y'\in A$. Without loss of generality, we suppose that for any  $z\in\ga'$, $|z-z_0|<|y-z_0|$. Since $\ga'\cap\ga=\0$, one of the liftings $\tga'$ of the arc $\ga'$ lies in some half-strip bounded by the lines $\tga$ and $\tga +2i\pi$ and a segment $[\tilde y, \tilde y+2i\pi]$ of the line $Re(w)=\rho_0$. Therefore, for any $w\in \tga'$, $|\te-\te_0-\pi-\la(\rho-\rho_0)|<M+\pi$. Applying similar argument to a subarc $\ga''\IN Q$ with the endpoint $z_0$ and the arc $\ga'$, we obtain the inequality for any $w\in \tga''$, $|\te-\te_0-2\pi-\la(\rho-\rho_0)|<M+2\pi$.\vse

Thus we arrive to the following theorem for FI(s)-systems of similarities on the plane.
\begin{thm}\label{parmat}
Let  $\eS=\{S_1, . . ., S_m\}$ be a FI(s)-system of  similarities in $\rr^2$. If one of the addresses of a point $p\in P$ is preperiodic,
 then  Jordan arcs $\ga\IN K$ with the end point $p$ have a slope parameter $\la(p)$, which is the same for all these arcs.
\end{thm}
\dok  Let $\bj\bar\bi$ be the preperiodic address of the point $p$.  Then the point $p'=S_\bj^{-1}(p)$ is a fixed point of $S_\bi$, so $\pim{p'}\ni \bar\bi$.

 There is $l\neq j_1$ such that $p\in K_{j_1}\cap K_l$. 
 Notice that there is a Jordan arc $\ga''\IN (K_l\mmm p)\in \bbc\mmm K_l$ whose endpoint is $p$. Therefore, $p'$ lies on the boundary of some component $A$ of $\bbc\mmm K$ which contains $S_\bj^{-1}(\ga'')$. By Proposition \ref{invarcs}, all the Jordan arcs $\ga\IN K$ with endpoint $p'$ have a slope parameter $\la$, which implies that all Jordan arcs $\ga\in K$ which have the endpoint at $p$ also have the slope parameter $\la$. \vse

\end{document}